\documentclass[12pt]{article}
\usepackage[intlimits]{amsmath}
\usepackage{amsthm,amssymb}
\usepackage{indentfirst}

\newtheorem{Theo}{Theorem}
\newtheorem{Lem}{Lemma}
\newtheorem{Prop}{Proposition}
\newtheorem{Cor}{Corollary}

\theoremstyle{definition}

\newtheorem{Def}{Definition}

\newtheorem{Rem}{Remark}

\theoremstyle{remark}

\def\cpict#1#2#3#4#5{\begin{figure}[h]
\dimen1=#1 \advance\dimen1 by -\hsize \divide\dimen1 by -2
\noindent\hskip\dimen1{\special{em:graph #3}} \vspace{#2}
\caption{#4}\label{#5}
\end{figure}}

\def\op{\operatorname}
\title{Boundary Value Problems on Manifolds with Fibered Boundary}
\author{A.~Savin, B.~Sternin}

\date{}
\begin{document}
\maketitle
\def\refname{References}

\section*{Introduction}
The modern theory of elliptic boundary value problems deals with
two types of problems. First, one has classical boundary value
problems that admit a realization as some bounded operators in
Sobolev spaces. Secondly, there are boundary value problems of
Atiyah--Patodi--Singer type that can be realized as Fredholm
operators in some \textit{subspaces} of the Sobolev spaces.
Moreover, the subspaces in question are the ranges of some
pseudodifferential projections acting in the Sobolev spaces. These
two classes of boundary value problems have a substantial
difference. Namely, the first class is defined \textit{only for
some} elliptic operators (defined on manifolds with boundary) and
the corresponding obstruction in the stable theory is known as the
Atiyah--Bott condition \cite{AtBo1}. The second class of boundary
value problems does not have this restriction: Fredholm boundary
value problems of the above described type can be defined for an
arbitrary elliptic operator. On the other hand, this type of
boundary conditions imposes an essential restriction on the
right-hand sides of the boundary value problem.  Namely, it is
supposed that the right-hand sides are taken from a subspace of
the Sobolev space of a possibly infinite codimension. The
following question naturally arises: is it possible to construct
an elliptic theory, which is a deformation of the two theories
such that on the one end it reduces to the theory of classical
boundary value problems, while at the other end it gives the
Atiyah--Patodi--Singer problems? In other words, the problem is to
construct a series of intermediate theories of elliptic boundary
value problems that would contain as particular (and in a sense
polar) examples classical boundary value problems and the problems
with projections. The present paper is devoted to this problem.

Clearly, the construction of such a theory requires some additional
assumptions on the geometry of the boundary of the manifold, where the
boundary value problems are considered. In the present paper we assume
that the boundary is a
\textit{fibration} over a compact base with a compact fiber.
In this setting, we establish all the analytical results we need
(finiteness theorem, in particular). Moreover, we compute the obstruction
to the existence of elliptic problems of this type. The obstruction turns out
to be an element of  $K^1$-group of the cotangent bundle of the base. This
explains the nature of the topological Atiyah--Bott obstruction as well
as the absence of obstruction for Atiyah--Patodi--Singer problems.
Indeed, in the special case, when the base coincides
with the boundary, one obtains classical boundary value problems.
On the contrary, when the base of the fibration is a point, we
obtain Atiyah--Patodi--Singer problems.

It should be mentioned that the importance and the interest of
this class of boundary value problems on manifolds with fibered
boundary is clear even in the case, when the boundary is a covering, i.e. a
fibration with a discrete fiber. For a covering, the class of
boundary value problems under consideration includes a number of
\textit{nonlocal} boundary value problems (see \cite{SaSt10}).
Rather surprisingly, the boundary conditions in the intermediate
theories on manifolds with fibered boundary are defined by
operators with discontinuous symbols.

Let us briefly describe the contents of the paper. We construct an
algebra of operators with discontinuous symbols on a fibration in
the first three sections. In particular, we establish the
composition formula and prove the Fredholm criterion. Boundary
value problems on manifolds with fibered boundary are defined in
Section \ref{s5}. We prove the finiteness theorem and in Section
\ref{s6}  compute the topological obstruction to pose a Fredholm
boundary value problem for an elliptic operator. Examples of
boundary value problems for the Hirzebruch operator are presented.

We are grateful to Professor G.~Rozenblioum of Chalmers University,
G\"oteborg, Sweden for numerous helpful discussions and
Dr.~V.~Nazaikinskii for help on a number of topics.
The results were
announced at the Conference "Workshop in Partial Differential
Equations", Potsdam, Germany,  November 12--16, 2001. We would like
to thank the organizers of this conference for their hospitality
and support. The work was partially supported by RFBF grants
00-01-00161, 01-01-06013, 99-01-01254. The paper was partially written
at Chalmers University of Technology and supported by a grant from the
Swedish Royal Academy of Sciences.

\section{Symbolic algebra}

Let $\pi\,:\,M\rightarrow X$ be a locally trivial fiber bundle of
compact smooth manifolds $M$ and $X$ with a typical fiber $Y$.
Manifolds $M,Y,X$ are assumed to be closed.
The fiber $\pi^{-1}(x)$ over $x\in X$ is denoted by $Y_x$.

On the total space of the fibration,  we consider special
coordinates described as follows.  For a domain $U\subset X$ with
coordinates $x_1,\ldots,x_n,\,\, n=\op{dim}X$ and a trivializing
mapping $\alpha_U$
\begin{align*}
  \pi^{-1}(U)&{\phantom{2}}\overset{\alpha_U}{\longrightarrow} \,\,U\times Y \\
  \quad \pi  &\searrow \quad  \swarrow \\ & {\phantom{11}}U
\end{align*}
we consider some coordinates $y_1,\ldots,y_m, \,\, m=\op{dim}Y$ in
a domain $V\subset Y$  of the fiber. Denote the dual coordinates
in the cotangent space by $(\xi_1,\ldots,
\xi_n,\eta_1,\ldots,\eta_m)$.

\begin{Def}\label{d1}
\textit{Principal symbol} of order zero is a function
\[
a_M\in C^\infty (T^*M\setminus \pi^*T^*X),
\]
which is homogeneous in covariables of order zero and is
directionally smooth as covariables  approach the plane
$\pi^*T^*X\subset T^*M$. More precisely, in local coordinates we
suppose that the function
\begin{equation}\label{*}
  a(x,y,\xi,t\eta), \quad (x,y,\xi,t,\eta)\in U\times
  V\times(\mathbb{R}^n\setminus\{0\})\times\mathbb{R}_+\times(\mathbb{R}^m\setminus\{0\})
\end{equation}
extends smoothly up to $t=0$.
\end{Def}
The following notation will be used for the directional limit
\[
\widetilde{a}_M(x,y,\xi,\eta)=\underset{t\to 0}{\op{lim}}\,\,
a_M(x,y,\xi,t\eta)
\]
of the principal symbol at the plane $\pi^*T^*X$.
This function is defined if both $|\xi|\neq 0$ and $|\eta|\neq 0$.
It is homogeneous of order zero with respect to $\xi$ and $\eta$.

It is clear that the restriction of the principal symbol to the
cosphere bundle  $S^*M$ can have discontinuities. However, our
definition is equivalent to the requirement that the restriction
of the principal symbol to $S^*M$ extends to a smooth function on
a compact manifold $\overline{S^*M\setminus
\pi^*S^*X}$ with boundary. This space is a compactification of
$S^*M\setminus \pi^*S^*X $ obtained by attaching to it sequences
$(x_i,y_i,\xi_i,\eta_i)$ converging to a point in $\pi^*S^*X$ such that the
quotients $\eta/|\eta|$ converge as well. One can show that this
manifold is diffeomorphic to $S^*M\setminus U_{\pi^*S^*X}$, where
$U_{\pi^*S^*X}$ is an open tubular neighborhood of $\pi^*S^*X$.

\begin{Def}\label{d2}
\textit{Operator symbol} of order zero is a function
\[
a_X\in C^\infty(T^*X\setminus \{0\},\Psi^0(Y))
\]
that at a point $(x,\xi)\in T^*X\setminus 0$ is a classical
pseudodifferential operator of order zero in the fiber $Y_x$. The
function $a_X$ is assumed to be homogeneous in covariable of order
zero and  smoothly depending on $(x,\xi)$. More precisely,
$a_X(x,\xi)$ is a smooth function with values in a Fr\'echet
space.
\end{Def}

\begin{Rem}\label{r1}
Let us recall the definition of the Fr\'echet structure on the
space of classical pseudodifferential operators, see
\cite{AtSi4,Schu4}. The set of seminorms is defined as follows.
Fix on $Y$ a quantization mapping $\sigma\mapsto
\widehat{\sigma}$, defined for smooth homogeneous symbols $\sigma$. Then a
classical $\psi DO$ $A$ of order $d$ on $Y$ is equivalent to a sum
\[
A\sim\sum_{j\geq 0}\widehat{a}_{d-j}, \quad a_{d-j}\in C^\infty
(S^*Y).
\]
Here the symbols $a_k$ are homogeneous of degree $k$. While the
equivalence means that
\[
A-\sum_{N\geq j\geq 0}\widehat{a}_{d-j}\in\Psi^{d-N-1}(Y).
\]
In these terms, the following two systems of seminorms on the
space $\Psi^d(Y)$:

1) $\|A\|_{\alpha,j}=\|a_{d-j}\|_\alpha$, where $\alpha$ runs over
a countable set of seminorms on $C^\infty (S^*Y)$;

2) $\|A\|_{s,N}=\|A-\sum_{0\leq j\leq
N-1}\widehat{a}_{d-j}\|_{H^s(Y)\rightarrow H^{s-(d-N)}(Y)}$. This
is an operator norm;

define the Fr\'echet structure.
\end{Rem}

\begin{Def}\label{d3}
Symbols $a_M$ and $a_X$ are said to be \textit{compatible} if
\begin{equation}\label{88}
\sigma(a_X)=\widetilde{a}_M,
\end{equation}
i.e. the principal symbol of the operator symbol $a_X$ is equal to
the limiting value $\widetilde{a}_M$ of the principal symbol at
the horizontal bundle $\pi^*T^*X$.
\end{Def}

\begin{Def}\label{d4}
By $\Sigma (M,\pi)$ denote the algebra of compatible pairs
$(a_M,a_X)$ with a componentwise product. An element
$\sigma\in\Sigma(M,\pi)$ is called a \textit{symbol} on $M$.
\end{Def}

\noindent {\bf Example 1.} Let $B\in C^\infty(X,\Psi^0(Y))$ be a
smooth family of pseudodifferential  operators $B_x$ of order zero
in the fibers:
  \[
B_x\,:\,C^\infty(Y_x)\rightarrow C^\infty(Y_x).
  \]
Then the pair
  \[
(\sigma(B),B)
  \]
is a symbol on $M$. Indeed, the principal symbol
$\sigma(B)(x,y,\eta)$ is smooth for $|\eta|\neq 0$, while the
rescaling as in~\eqref{*} does not change the symbol:
\[
\sigma(B)(x,y,t\eta)=\sigma (B)(x,y,\eta),
\]
since $\sigma (B)$ is homogeneous of order zero. Thus, in this
case $\widetilde{\sigma(B)}=\sigma(B)$ and the compatibility
condition~\eqref{88} is clearly satisfied.

\vspace{0.4cm}\noindent {\bf Example 2.} Let us now choose a
smooth symbol $a\in C^\infty (T^*M\setminus
  \{0\})$. Then the pair
  \[
(a, a|_{\pi^* T^* X}),
  \]
where operator symbol $a\big|_{\pi^*T^*X}=a(x,y,\xi,0)$ acts as a
multiplication operator, is a symbol in the sense of Definition \ref{d4}.

Denote by $\Sigma_0\subset\Sigma (M,\pi)$ the subalgebra
multiplicatively generated by the symbols from Examples 1 and 2.

\begin{Prop}\label{pr1}
$\Sigma_0$ is dense in $\Sigma(M,\pi)$, where $\Sigma (M,\pi)$ is
equipped with the Fr\'echet structure as a subspace:
\[
\Sigma(M,\pi)\subset C^\infty (S^*M\setminus U_{\pi^*S^*X})\oplus
C^\infty (S^*X,\Psi^0 (Y)).
\]
\end{Prop}

\begin{proof}

1) Let us show that the subalgebra of
  principal symbols corresponding to $\Sigma_0$ is dense.

Consider a principal symbol $\sigma (x,y,\xi,\eta)$. By a
partition of unity argument, it suffices to approximate
\footnote{Here and in what follows we use the term approximation
in the sense of Fr\'echet topologies on the corresponding spaces.}
$\sigma$ over a neighborhood of an arbitrary point $(x_0,y_0)\in
M$, where we can use local coordinates. In the domain
$|\eta|>2\varepsilon\sqrt{\xi^2 + \eta^2}$, where $\sigma$ is
smooth, it is equal to a smooth symbol
  \[
\chi\biggl(\frac{|\eta|}{\sqrt{\xi^2 +\eta^2}}\biggr)\,\sigma
(x,y,\xi,\eta),
  \]
where the chopping function $\chi(t)$ is zero for $t<\varepsilon$
and is $1$ for $t>2\varepsilon$. Thus, the main problem is to
approximate the symbol for small $\eta$.

To this end, consider the function $\sigma(x,y,\xi', t\eta')$.
It is smooth for
\begin{equation}\label{+}
|\xi'|=|\eta'|=1  \quad\text{and}\,\,\, t\in [0,100].
\end{equation}
Hence, it can be approximated by a sum of products of functions
depending on two subsets of variables $(x,\xi',t)$ and $(y,\xi')$:
\[
\sigma(x,y,\xi', t\eta')\sim\sum_{j=1}^N a_j
(x,\xi',t)b_j(y,\eta'),
\]
where the terms are smooth functions for the parameter values as
in~\eqref{+}. Thus, the original symbol is approximated for small
$\eta$ by the expression

\begin{equation}\label{e1}
\sigma(x,y,\xi, \eta)\sim\sum_{j=1}^N
a_j\biggl(x,\frac{\xi}{|\xi|},
\frac{|\eta|}{|\xi|}\biggr)b_j\biggl(y, \frac{\eta}{|\eta|}\biggr)
\end{equation}
for $|\eta|\leq 100|\xi|$. Here the terms $b_j$ have the desired
form as in Example 1 and they belong to the subalgebra $\Sigma_0$.
However, the functions $a_j$ are not smooth at $\eta=0$ (they
contain terms $|\eta|$). Let us eliminate this singularity.

To this end, consider the Taylor expansion of $a_j(x,\xi',t)$ in
$t$ at zero. We obtain
\begin{equation}\label{e2}
a_j\biggl(x,\frac{\xi}{|\xi|},
\frac{|\eta|}{|\xi|}\biggr)=\sum^{N'}_{k=0}\frac{\partial
^k}{\partial t^k}a_j\biggl(x,
\frac{\xi}{|\xi|},0\biggr)\biggl(\frac{|\eta|}{|\xi|}\biggr)^k\frac{1}{k!}+O\biggl(
\frac{|\eta|^{N'+1}}{|\xi|^{N'+1}}\biggr).
\end{equation}
For $k$ even the corresponding term $|\eta|^k/|\xi|^k$ is smooth
for $|\xi|\neq 0$. However, for $k$ odd the quotient
$|\eta|/|\xi|$ is not  smooth.  We represent it in the form
\begin{equation}\label{e3}
\frac{|\eta|}{|\xi|}=\sum_i\frac{\varphi_i
\biggl(\dfrac{\eta}{|\eta|}\biggr)|\eta|}{\eta_i}
\frac{\eta_i}{|\xi|},
  \end{equation}
where $\varphi_i$ denotes a partition of unity on the sphere
$|\eta|=1$ such that $\eta_i$ does not vanish on the support of
$\varphi_i$. In this sum, the factors on the left correspond to
symbols from Example 1, while the quotient $\eta_i/|\xi|$ is
smooth and corresponds to a symbol from Example 2.

Now the approximation of the symbol $\sigma$ is defined as
follows. First, we approximate $\sigma$ for $|\eta|<|\xi|$ by the
formula~\eqref{e1}. Then the Taylor expansion in~\eqref{e2} is
written with the decomposition of nonsmooth factors as
in~\eqref{e3}. The error term in the Taylor expansion shows that
for $|\eta|<4\varepsilon |\xi|$ (for $\varepsilon$ sufficiently
small), the approximation is good. Denote the corresponding
approximating element by
\[
\sum_\alpha a_\alpha(x,\xi,\eta)b_\alpha(y,\eta)
\]
($a_\alpha$ and $b_\alpha$ are homogeneous in $(\xi,\eta)$ and
$\eta$, respectively).  Then the expression
\[
\chi\biggl(\frac{|\eta|}{\sqrt{\xi^2
+\eta^2}}\biggr)\sigma(x,y,\xi,\eta)+\sum_\alpha\biggl(\bigg[1-
\chi\biggl(\frac{|\eta|}{\sqrt{\xi^2
+\eta^2}}\biggr)\bigg]a_\alpha(x,\xi,\eta)\biggr)b_\alpha(y,\eta)
\]
approximates the initial principal symbol with a small error as
desired. This approximation can now be globalized in $x$ and $y$
by a partition of unity.

2) Let us show that a compatible pair $(p_M,p_X)$ of a principal
symbol and an operator symbol can be approximated by elements from
the subalgebra $\Sigma_0$.

It follows from the previous part of the proof  that we need to
prove this for the trivial principal symbol $p_M=0$ only.

So, let us consider a symbol $(0,p_X)$ defined by an
operator-valued function $p_X(x,\xi)\in\Psi^{-1}(Y_x)$.

Similar to the above, let us approximate $p_X$ in a neighborhood
of a point $(x_0,\xi_0)\in S^*X$. There $p_X$ can be approximated
by its Taylor expansion in $x$ and $\xi$. We write this for
brevity as
  \[
p_X(x,\xi)=\sum_{1\leq j\leq N'} p_{j,x_0,\xi_0}(x,\xi)B_j
+O(|x-x_0|^N+|\xi-\xi_0|^N),
  \]
where $B_j$ are pseudodifferential operators in the fibers and
$p_{j,x_0,\xi_0}(x,\xi)$ are smooth (scalar) functions.

Taking a cover of $S^*M$ with sufficiently small charts $U_\alpha$
and using the corresponding partition of unity $\varphi_\alpha$,
we obtain an approximation
\[
\sum_{\alpha,j}\varphi_\alpha p_{j,x_\alpha,\xi_\alpha}(x,\xi)B_j
\]
for $p_X$.  If the charts $U_\alpha$ are chosen good enough
(e.\,g., their diameters are less than $\varepsilon$ and every
point belongs to at most $l+1$ charts, where $l$ is the dimension
of $S^*X$) then this expression is close to $p_X$ in the Fr\'echet
space $\Sigma (M,\pi)$. To end the proof, it suffices to show that
a term of the form $\varphi p(x,\xi)B$ is equal to a composition
of a smooth symbol and a symbol of a family of pseudodifferential
operators. Indeed, take $\chi\in C^\infty (X)$ such that $\pi^*
\chi\equiv 1$ in a neighborhood of the support of $\varphi$. Then
\[
\varphi p(x,\xi)B=[\varphi p(x,\xi)]\chi B.
\]
Here $\chi B$ is a family of $\psi DO$ in the fibers, while
$\varphi p(x,\xi)$ is the value for $\eta=0$ of the smooth symbol
\[
\chi (|\eta|)\varphi(x,\xi)p(x,\xi)
\]
where $\chi(t)$ is zero for $|t|>2\varepsilon$ and $1$ for
$|t|<\varepsilon$.  This completes the proof of the proposition.

\end{proof}

\section{Operators on fibered manifolds}\label{s2}

The aim of this section is to show that a symbol $\sigma
\in\Sigma(M,\pi)$ defines an operator
\[
\widehat{\sigma}\,:\,C^\infty (M)\longrightarrow C^\infty (M)
\]
acting on $C^\infty$ functions on $M$. Let us construct it. For a
decomposition $\sigma=(a_M,a_X)$ this operator is defined as
\begin{equation}
\label{kaa}
\widehat{\sigma}u=(\widehat{a}_M+\widehat{a_X-\widehat{\widetilde{a}}_M})u,
\quad u\in C^\infty(M),
\end{equation}
where the first component $\widehat{a}_M$ denotes the usual
quantization of the principal symbol, while the second corresponds
to quantization of operator-valued symbols (note that there is a
correction term $\widehat{\widetilde{a}}_M$ to the operator
symbol. This will be described later on).

More precisely, $\widehat{a}_M$ is first defined over charts
$U_\alpha\subset X$ by the formula
\begin{align}\label{A}
A_\alpha u=\sum_j \frac{1}{(2\pi)^{m+n}}\psi_j(y)\int_{T^*_{x,y}M}
e^{i(x\xi+y\eta)}a_M(x,y,\xi,\eta)\widehat{(u\chi_j)}(\xi,\eta)\,d\xi
d\eta,
\end{align}
where $\chi_j(y)$ is a partition of unity subordinate to an atlas
in the fiber, $\psi_j$ is equal to $1$ in a neighborhood of the
support of $\chi_j$ and vanishes far from it. Finally, the hat
denotes the Fourier transform in $x$ and $y$.

The operator symbol $\widehat{\widetilde{a}}_M$ is defined
along the same lines
\[
\widehat{\widetilde{a}}_M(x,\xi)v=\frac{1}{(2\pi)^m}\sum_j
\psi_j(y)\int_{T^*_y Y_x}
e^{iy\eta}\widetilde{a}_M(x,y,\xi,\eta)\widehat{(v\chi_j)}(\eta)\,d\eta,
\quad v\in C^\infty(Y_x),
\]
quantizing only in the fiberwise variables. This definition uses
the same data $\chi_j,\psi_j$ as (\ref{A}).

In a chart $U_\alpha\subset X$, the operator
$\widehat{a_X-\widehat{\widetilde{a}}_M}$ with operator-valued
symbol $a_X-\widehat{\widetilde{a}}_M$ is defined as
\begin{equation}\label{b}
 B_\alpha u= \frac{1}{(2\pi)^n}\int_{T^*_{x}X}
e^{ix\xi}\big[a_X(x,\xi)-\widehat{\widetilde{a}}_M(x,\xi)\big]
\widehat{u}(\xi)\,d\xi.
\end{equation}
Here $\widehat{u}(\xi)$ denotes the Fourier transform in $x$ of $u(x,y)$.

Globally  on $M$ the operator $\widehat{\sigma}$ is defined using
the local expressions $A_\alpha, B_\alpha$:
\begin{equation}\label{*0}
\widehat{\sigma}u=\sum_\alpha \psi'_\alpha(A_\alpha+B_
\alpha)(\varphi'_\alpha u).
\end{equation}
by a partition of unity
$\chi'_\alpha$ on $X$ subordinate to the atlas $\{U_\alpha\}$
and functions $\psi'_\alpha$ supported in $U\alpha$ with the
property $\psi'_\alpha\chi'_\alpha=\chi'_\alpha$.

\begin{Rem}
If $\sigma=(a,a|_{\pi^*T^*X})$ corresponds to a smooth symbol as
in Example 2, then $\widehat{\sigma}$ is just the usual $\psi DO$
on $M$ with symbol $a$. On the other hand, for a family $B$ as in
Example 1 consider the symbol $\sigma=(\sigma(B),B)$. Then one can
prove that $\widehat{\sigma}=B$.
\end{Rem}

The quantization formula (\ref{kaa}) resembles the formula for operators
with discontinuous symbols on the circle, see \cite{PlRo2}.

\begin{Theo}\label{t1}
$\widehat{\sigma}$ extends to a bounded operator in Sobolev spaces
$H^s(M)$, $s\in \Bbb{R}$
\[
\widehat{\sigma}\,:\,H^s(M)\rightarrow H^s(M).
\]
\end{Theo}

\begin{proof}
1) The continuity in $H^s(M)$ of the first component~\eqref{A}
corresponding to the principal symbol was proved already in
\cite{KoNi1}. It was shown that no regularity of the symbol in
covariables is required to obtain continuity  in Sobolev spaces
(for all $s$). Moreover, the norm of an operator with symbol
$a(x',\xi')$ in $H^s(M)$ is bounded by the maximum of a finite
number of derivatives in the geometric variables:
\begin{equation}\label{Om}
 \underset{|\alpha|\leq N}{\op{max}}\,\,\, \underset{(x',\xi')\in T^* M}{
 \op{sup}}\big|\partial^\alpha_{x'}a(x',\xi')\big|,
\end{equation}
up to a constant depending only on $s$ (the order $N$ of the derivatives
also depends on $s$).

2) The continuity of operator~\eqref{b} with operator-valued
symbol in $L^2$ follows from the paper \cite{Luk1}. In this case
the $L^2$-norm is also bounded by estimates of the
form~\eqref{Om}, where $a$ is replaced by the corresponding
operator-valued symbol and the absolute value is replaces by the operator
norm.

3) Let us prove that the term ~\eqref{b} defines a bounded
operator in the Sobolev spaces. By interpolation it suffices to
prove the boundedness for $s=\pm 2m,\,\, m\in \mathbb{N}$.

If $s=2m$ then according to Theorem \ref{bounded} of the
Appendix an operator
\[
\widehat{q}_X\,:\,C^ \infty (M)\rightarrow C^\infty (M)
\]
with operator-valued symbol $q_X(x,\xi)\in\Psi^{-1}(Y_x)$
$$
q_X=a_X-\widehat{\widetilde{a}}_M
$$
acts continuously in $H^{2m}(M)$ if the composition
\begin{equation}\label{=}
  (\triangle_Y+\xi^2)^m{q}_X(x,\xi)
  (\triangle_Y+\xi^2)^{-m}\,:\, C^\infty(Y)\rightarrow
  C^\infty (Y)
\end{equation}
is uniformly bounded in $L^2 (Y)$, where $\triangle_Y$ is a positive
Laplacian on $Y$.

The corresponding estimate
\[
\|(\xi^2+\triangle_Y)^m
q_X(x,\xi)(\xi^2+\triangle_Y)^{-m}\|_{L^2(Y)\rightarrow
L^2(Y)}\leq C
\]
can be proved using a decomposition of the commutator
\begin{equation}
\left[(\triangle_Y+\xi^2)^m,q_X\right]=\sum_{-2\le j \le 2m-2}
a_j(x,\xi),
\end{equation}
where the component $a_j$ has order $j$ and its norm in the spaces
$$
a_j(x,\xi):H^s(Y)\longrightarrow H^{s-j}(Y)
$$
is bounded by $(1+|\xi|)^{2m-2-j}$. This decomposition can be proved by
induction.

Using this decomposition, we can estimate each of the terms:
$$
\|a_j(x,\xi)(\xi^2+\triangle_Y)^{-m}\|_{L^2(Y)\rightarrow
L^2(Y)}\le\|a_j\|_{H^{j}(Y)\to H^0(Y)}\cdot
\|(\triangle_Y+\xi^2)^{-m}\|_{H^0(Y)\to H^{j}(Y)}.
$$
Thus, we obtain
$$
\|a_j(x,\xi)(\xi^2+\triangle_Y)^{-m}\|_{L^2(Y)\rightarrow L^2(Y)}\le
(1+|\xi|)^{2m-2-j}\cdot {\rm max}
\left((1+|\xi|)^{-2m+j},(1+|\xi|)^{-2m}\right).
$$
This expression is clearly uniformly bounded.

The remaining case $s=-2m$ can be considered similarly.

\end{proof}

\section{Operator algebra}\label{s3}

\begin{Theo}\label{t2}
The following composition formula is valid
\[
\widehat{\sigma}_1\widehat{\sigma}_2=\widehat{\sigma_1\sigma_2}+K,
\]
where $\sigma_1,\sigma_2\in\Sigma(M,\pi)$ are symbols, while the
error term $K$ is a compact operator in all spaces $H^s(M)$.
Moreover, if $\sigma_1\in\Sigma_0$ then $K$ has order $-1$
in the scale.
\end{Theo}

The proof of this theorem is done in the Appendix.

\begin{Def}
Denote by $\Psi^0(M,\pi)$ the space of operators of the form
\[
\widehat{\sigma}+K\,:\,C^\infty (M)\longrightarrow C^\infty (M),
\]
where operator  $K$ extends to a compact operator in the scale
$H^s(M)$.

The composition formula enables us to prove the following result.
\end{Def}

\begin{Theo}\label{t3}

1) $\Psi^0(M,\pi)$ is an algebra.

2) The subalgebra $\Psi_0\subset \Psi^0(M,\pi)$ generated by $\psi
DO$'s on $M$, families of $\psi DO$'s in the fibers and compact
operators is dense in $\Psi^0(M,\pi)$ with respect to operator
norm in $H^s(M)$.

3) The norm, modulo compact operators, is equal to
\begin{equation}\label{*0*}
\underset{k\in
K(L^2(M))}{\op{inf}}\|\widehat{\sigma}+k\|=\op{max}\biggl(\underset{(x,y,\xi,
\eta)\in
S^*M\setminus\pi^*S^*X}{\op{sup}}\big|a_M(x,y,\xi,\eta)\big|,\,
\underset{(x,\xi)\in S^* X}{\op{max}}\|a_X(x,\xi)\|\biggr),
\end{equation}
where $K(L^2(M))$ denotes the ideal of compact operators.
\end{Theo}

\begin{proof}
1) This straightforwardly follows from Theorem \ref{t2}.

2) Suppose that $\widehat{\sigma}\in \Psi^0(M,\pi)$. According to
Proposition \ref{pr1} its symbol $\sigma$ can be
approximated by a symbol $\sigma_\varepsilon \in \Sigma{}_0$. Then
by Theorem \ref{t2} and Remark 2 we have
$\widehat{\sigma}_\varepsilon\in\Psi_0$. On the other hand, the
difference $\widehat{\sigma}-\widehat{\sigma}_\varepsilon=
\widehat{\sigma-\sigma_\varepsilon}$
has a small symbol. Thus, its norm is small. This proves that
$\Psi_0$ is dense in $\Psi^0(M,\pi)$.

3) By virtue  of the second part of the theorem, it suffices to
prove the estimate for $\widehat{\sigma}\in\Psi_0$ and
$\sigma=(a_M,a_X)$.

Let us first prove the estimate from below. Suppose that
$(x_0,y_0,\xi_0,\eta_0)\in T^*M\setminus \pi^*T^*X$. Then we
choose a sequence of smooth functions
\[
u_n(x,y,t)=e^{it(x\xi_0+y\eta_0)}\chi_n(x,y),
\]
where $\|\chi_n\|_{L^2(M)}=1$ and $\chi_n$ is supported in a ball
of radius $ 1/n$ around $(x_0,y_0)$.

Since $\widehat{\sigma}\in\Psi_0$, it has the form (modulo a
compact operator)
\begin{equation}\label{w}
\widehat{\sigma}=\sum_\alpha \prod_\beta
A_{\alpha\beta}B_{\alpha\beta},
\end{equation}
where $A_{\alpha\beta}$ are $\psi DO$'s on $M$ and
$B_{\alpha\beta}$ are families of $\psi DO$'s in the fibers as
in Example 1.

Then as $t$ tends to infinity we obtain
\[
A_{\alpha\beta}u_n=\sigma(A_{\alpha\beta})(x,y,\xi_0,\eta_0)u_n+O\biggl(
\frac 1t\biggr),
\]
\[
B_{\alpha\beta}u_n=\sigma(B_{\alpha\beta})(x,y,\eta_0)u_n+O\biggl(
\frac 1t\biggr),
\]
according to the H\"ormander's definition of pseudodifferential
operators.

Hence, as $t\to \infty$ we obtain
\[
\widehat{\sigma}
u_n=a_M(x,y,\xi_0,\eta_0)u_n+O
\biggl(\frac 1t\biggr).
\]
On the other hand,  $u_n$ weakly converges to zero as $t\to
\infty$. Thus, for a compact operator $K$ we have
$K{u_n}\rightarrow 0$. Therefore,
\[
(\widehat{\sigma}+K)u_n=a_M(x_0,y_0,\xi_0,\eta_0)u_n+u_\varepsilon+O
\biggl(\frac 1t\biggr)+Ku_n,
\]
where $u_\varepsilon
=(a_M(x,y,\xi,\eta)-a_M(x_0,y_0,\xi_0,\eta_0))u_n$, and
for $n$ large we have $\|u_\varepsilon\|<\varepsilon$. For $t$ large
this yields
\[
\|(\widehat{\sigma}+K)u_n\|_{L^2(M)}\geq
|a_M(x_0,y_0,\xi_0,\eta_0)|\cdot\|u_n\|_{L^2(M)}-\varepsilon.
\]
Since $\varepsilon$ can be chosen arbitrarily small and
$(x_0,y_0,\xi_0,\eta_0)$ and $K$ are also arbitrary, we have the
desired estimate of the norm modulo compact operators
\[
\underset{K}{\op{inf}}\|\widehat{\sigma}+K\|\geq
\underset{(x,y,\xi,\eta)\in T^*M\setminus
\pi^*T^*X}{\op{sup}}|a_M(x,y,\xi,\eta)|.
\]

Let us  prove the second part of the lower  estimate that deals
with the operator symbol. For $(x_0,\xi_0)\in S^*X$ and $v\in
C^\infty(Y_{x_0})$ consider the sequence
\[
u_n(x,y)=e^{ix\xi_0 t} v(y)\chi_n(x),
\]
where $\chi_n(x)$, $\|\chi_n\|_{L^2(X)}=1$ is a function on
the base $X$ supported in a ball of radius $1/n$ around $x_0$.

For $\widehat{\sigma}$ as in~\eqref{w}, we obtain
\[
A_{\alpha\beta}u_n=\sigma(A)(x,y,\xi_0,0)u_n+O\biggl(\frac
1t\biggr),
\]
and
\[
B_{\alpha\beta}u_n=B_{\alpha\beta}u_n.
\]
By an argument similar to the previous part of the proof, we
obtain
\[
\underset{k\in K(L^2(M))}{\op{inf}}\|\widehat{\sigma}+k\|\geq
\underset{(x,\xi)\in
S^*X}{\op{max}}\|a_X(x,\xi)\|_{L^2(Y_x)\rightarrow L^2(Y_x)}.
\]

We now prove the estimate from above. This is done following
the standard scheme (see, e.\,g. \cite{ReSc1}).

Let us denote by $\|\sigma\|$ the norm of the symbol
$\sigma=(a_M,a_X)$, as defined by the right hand side of~\eqref{*0*}.
Fix a constant $C>\|\sigma\|$. Then $C^2-\sigma^*\sigma$ is self-adjoint
and positive symbol. Denote by $\sigma_0$ its positive square root
\[
C^2-\sigma^*\sigma=\sigma^2_0.
\]
Let us now approximate $\sigma_0$ by a symbol
$\sigma_\varepsilon\in\Sigma_0$. Thus, according to the
composition formula, we have for the corresponding operators
\[
C^2-\widehat{\sigma}^*\widehat{\sigma}=\widehat{\sigma}_{\varepsilon}^*
\widehat{\sigma}_\varepsilon+K_{-1}+N_\varepsilon,
\]
where $K_{-1}$ is a compact operator of order $-1$
in the Sobolev scale and $N_\varepsilon$ has norm less than
$\varepsilon$. This
formula gives the following estimate
\begin{equation}\label{AA}
  \|\widehat{\sigma}u\|^2_{L^2(M)}\leq
  C^2\|u\|^2_{L^2(M)}+\|K_{-1} u\|_{L^2(M)}\cdot\|u\|_{L^2(M)}+\varepsilon
  \|u\|^2_{L^2(M)}.
\end{equation}
Denote by $R_\varepsilon$ the smoothing operator on $L^2(M)$ with
symbol
\[
\sigma(R_\varepsilon)(x,y,\xi,\eta)=
  \begin{cases}
    0, & \text{for}\,\, |\xi|^2+|\eta|^2>2/\varepsilon, \\
    1, & \text{for}\,\, |\xi|^2+|\eta|^2<1/\varepsilon.
  \end{cases}
\]
It  is a compact operator and the following estimates are valid
(see \cite{ReSc1})
\begin{equation}\label{B}
\|u-R_\varepsilon u\|_{L^2(M)}\leq \|u\|_{L^2(M)} \quad \text{and}
\quad \|u-R_\varepsilon u\|_{H^{-1}(M)}\leq
c(\varepsilon)\|u\|_{L^2(M)},
\end{equation}
where $c(\varepsilon)\to 0$ as $\varepsilon \to 0$.
Then~\eqref{AA} and~\eqref{B} give the desired estimate:
\[
\|\widehat{\sigma}(u-R_\varepsilon u)\|^2_{L^2(M)}\leq
C^2\|u\|^2_{L^2(M)}+C(\varepsilon)\|u\|^2_{L^2(M)}
\]
with $C(\varepsilon)\to 0$ as $\varepsilon\to 0.$
Hence, if we take  the compact operator
$k=\widehat{\sigma}R_\varepsilon$, this yields for the infinum of
the norm the following estimate
\[
\underset{k\in K(L^2(M))}{\op{inf}}\|\widehat{\sigma}+k\|\leq
\|\sigma\|
\]
as desired.

Theorem is proved.

\end{proof}

Denote by $\overline{\Psi}^0(M,\pi)$ the closure of
$\Psi^0(M,\pi)$ with regard for the operator norm in $L^2(M)$. The
next result describes the corresponding Calkin algebra.

\begin{Cor}
\[
\overline{\Psi}^0(M,\pi)/\mathcal{K}\simeq\overline{\Sigma}(M,\pi),
\]
where $\mathcal{K}$ denotes the ideal of compact operators and
$\overline{\Sigma}(M,\pi)$ is the completion of the symbolic
algebra ${\Sigma}(M,\pi)$ with respect to the norm $\|\cdot\|$. In
more detail, the completion is a subalgebra
\[
\overline{\Sigma}(M,\pi)\subset C(S^*M\setminus
U_{\pi^*S^*X})\oplus C(S^*X,\overline{\Psi}^0(Y)),
\]
consisting of compatible pairs $(a_M,a_X)$:
\[
\overline{\Sigma}(M,\pi)=\bigr\{a_M\oplus a_X \big|\quad
a_M\big|_{\partial U_{\pi^*S^*X}}=\sigma(a_X)\bigl\},
\]
where $U_{\pi^*S^*X}$ denotes an open tubular neighborhood of
$\pi^*S^*X$ in $S^*M$ and $\overline{\Psi}^0(Y)$ is the norm
closure of the algebra of pseudodifferential operators of order
zero in the fibers.
\end{Cor}

The proof of the corollary follows easily from the estimate modulo
compact operators.

The description of the Calkin algebra enables us to state the
Fredholm criterion.

\begin{Cor}
$\widehat{\sigma}\in \Psi^0(M,\pi)$ is a Fredholm operator in
Sobolev spaces if and only if $\sigma$ is invertible.
\end{Cor}

\begin{proof}
The ``if'' part follows from the composition formula: the
parametrix is given by $\widehat{\sigma^{-1}}$. The ``only if''
part is proved as follows. Let us assume at first that $s=0$.
Suppose that $\widehat{\sigma}$ is a Fredholm operator in
$L^2(M)$ with a left quasiinverse $A$:
\[
A\widehat{\sigma}=1+K_1.
\]
Then this gives an apriori estimate
\[
\|u\|\leq C \|\widehat{\sigma} u\|+\|K_1 u\|.
\]
If we substitute in this inequality the sequence $u_n(t)$, as in the
proof of Theorem \ref{t3}, we obtain, choosing $n$ big enough and letting
$t\to \infty$, that
\[
1\leq C_1 |\sigma_M (x_0,y_0,\xi_0,\eta_0)|
\]
and a similar estimate for the operator symbol
\[
1 \leq C_2 \|\sigma_X(x_0,\xi_0)\|
\]
(these estimates can be obtained choosing approximations of
$\widehat{\sigma}$ by elements of the subalgebra $\Psi_0$).
Thus, the symbol is monomorphic. Passing to the adjoint operator,
one proves the surjectivity.

Therefore, the symbol of a Fredholm operator is an isomorphism.
The ellipticity of Fredholm operators in $H^s(M),\,s\neq 0$ can be proved
along the same lines using the compositions
$\triangle_M^{s/2}\widehat{\sigma}\triangle_M^{-s/2}$.
\end{proof}

\begin{Rem}
A generalization of the composition formula and the finiteness
theorem to operators acting in the sections of some vector bundles
over $M$ is rather standard and is left to the reader, e.\,g. see
\cite{AtSi4}.
\end{Rem}

\section{Elliptic Theory in subspaces}\label{s4}
Let us consider operators
\begin{equation}\label{ww}
  D\,:\,\op{Im} P_1\longrightarrow \op{Im} P_2
\end{equation}
acting in subspaces defined by projections
\[
P_1\,:\,C^\infty (M,E)\longrightarrow C^\infty (M,E), \quad
P_2\,:\, C^\infty (M,F)\longrightarrow C^\infty (M,F),
\]
where $E,F$ are vector bundles over $M$.
We suppose that $P_{1,2}$ belong  to the algebra defined in
previous sections: $P_{1,2}\in \Psi^0(M,\pi)$ and operator $D$ is
a restriction of some  operator $\widetilde{D}$
\[
\widetilde{D}\,:\, C^\infty(M,E)\longrightarrow C^\infty (M,F),
\]
also from our algebra: $\widetilde{D}\in\Psi^0(M,\pi)$.

\begin{Theo}\label{t5}
Operator~\eqref{ww} defines a Fredholm operator in Sobolev spaces
if and only if the following two conditions are satisfied:

1) the principal symbol
$$
\sigma_M(D)\,:\, \op{Im}\sigma_M
(P_1)\rightarrow \op{Im}\sigma_M (P_2)
$$
is invertible over $S^*M\setminus\pi^*S^*X$;

2) the operator symbol
$$
\sigma_X(D)\,:\, \op{Im}\sigma_X
(P_1)\rightarrow \op{Im}\sigma_X (P_2)
$$
is invertible over $S^*X$.
\end{Theo}

\begin{Rem}
In more detail, these conditions require that the principal symbol
is a  vector bundle isomorphism and the operator symbol is an
invertible family of operators in subspaces defined by pseudodifferential
projections (see, e.\,g., \cite{SaSt1}). Let us also mention that the
symbol of a projection is a projection itself.
\end{Rem}

\begin{proof}
If the two symbols are invertible then the parametrix has the
form:
\[
D^{-1}= P_1\widehat{\sigma^{-1}}\,:\,\op{Im} P_2\longrightarrow
\op{Im} P_1,
\]
where $\sigma=(\sigma_M(D),\sigma_X(D))$ is the symbol of $D$. The
fact that the differences $D^{-1}D-1$ and $DD^{-1}-1$ are compact
operators follows from the composition formula.

Let us prove the ``only if'' part. If $D$ is a Fredholm operator,
then the direct sum
\[
\op{Im} P_1\oplus \op{Im}P_1^\bot\overset{D^*D\oplus
1}{\longrightarrow}\op{Im}P_1\oplus\op{Im}P_1^\bot
\]
also has the Fredholm property. Moreover, there is an obvious
isomorphism
\[
\op{Im}P_1\oplus\op{Im}P_1^\bot\simeq H^s(M,E).
\]
Hence, $D^*D\oplus 1$ is an elliptic operator and its symbol has
trivial kernel. The surjectivity of the symbol is proved along the
same lines using the composition $DD^*$.

\end{proof}

\begin{Rem}
Consider a subclass of elliptic operators in subspaces~\eqref{ww},
where both projections $P_1,P_2$ correspond to families of
pseudodifferential projections in the fibers.  Then for the
identity map $\pi=\op{id}\,:\,M\rightarrow M$ we obtain the usual
operators on vector bundle sections over $M$, while for the
collapsing map $\pi\,:\,M\rightarrow pt$ this construction gives
the class of operators acting in subspaces defined by
pseudodifferential projections (see \cite{SaSt1}).
\end{Rem}

\section{Boundary value problems on manifolds with fibered
boundary}\label{s5}

{\bf 1. Main definitions.} Consider $M$ a compact smooth manifold
with boundary denoted by $\partial M$. Assume that $\partial M$ is the total
space of a locally trivial fiber bundle
\[
\pi\,:\,\partial M \longrightarrow X
\]
with a compact base $X$ and a compact fiber $Y$ as in the previous
sections.

On $M$ we consider an elliptic differential operator
\[
D\,:\,C^\infty (M,E)\longrightarrow C^\infty (M,F)
\]
of order $d$. To define the boundary conditions for $D$, we
introduce the operator
\[
j\,:\, C^\infty (M,E)\longrightarrow C^\infty
(M,E\big|^{d}_{\partial M}),
\]
that maps a function to its jet of order $d$ in the normal
direction to the boundary:
\[
ju=\biggl(u\big|_{\partial M}, -i\frac{\partial}{\partial
t}u\bigg|_{\partial M},\ldots, \biggl(-i \frac{\partial}{\partial
t}\biggr)^{d-1}u\bigg|_{\partial M}\biggr).
\]
Here $t\geq 0$ denotes the normal coordinate near the boundary. We
consider boundary value problems of the following form:
\begin{equation}\label{R}
  \begin{cases}
    Du=f, & u\in C^\infty(M,E),\qquad\quad\,\,  f\in C^\infty(M,F), \\
    Bju=g\in \op{Im}P, & \op{Im}P\subset C^\infty (\partial
M,G), \quad G\in \op{Vect}{(\partial M)},
  \end{cases}
  \end{equation}
where the subspace $\op{Im}P$  in the space of sections of a
vector bundle $G$ is defined by a family of pseudodifferential
projections over $X$ acting on functions in the fibers of $\pi$.

If $D$ is a first-order operator, then the boundary condition $B$
is assumed to be an element of the algebra $\Psi^0(\partial
M,\pi)$.  For operators of higher order the boundary condition
is more complicated. In this case it has $d$ components
\[
B\,:\,C^\infty(\partial M,E\big|^d_{\partial M})\rightarrow
C^\infty (\partial M,G)
\]
and is defined as a composition
\[
B=(B_0,B_1\triangle^{-1/2}_{\partial
M},\ldots,B_{d-1}\triangle^{-\frac{d-1}{2}}_{\partial M}),
\]
where $\triangle_{\partial M}$ is a positive Laplacian on
$\partial M$ and the components $B_j$ belong to $\Psi_0(\partial
M,\pi)$.

\begin{Rem}
This construction reduces to some well-known classes of boundary
value problems for special types of projections $\pi$:

1) $\pi=\op{id}\,:\,\partial M\rightarrow \partial M$. In this
case we obtain classical boundary value problems;

2) $\pi\,:\,\partial M\rightarrow pt$. This gives boundary value
problems in subspaces (see \cite{ScSS18});

3) $\pi\,:\, \partial M\rightarrow X$ and $\pi$ is a covering.
This gives a class of nonlocal boundary value problems studied in
\cite{SaSt10}.

\end{Rem}

{\bf 2. Finiteness theorem.}
The ellipticity condition of a boundary value problem for $D$ is
formulated in terms of a special vector bundle
\[
L_+(D)\in\op{Vect}(S^*\partial M)
\]
defined as a subbundle in the pull-back of $E\big|^d_{\partial M}$
to $S^*\partial M$ and generated by the Cauchy  data at $t=0$ of
functions $u(t)$ satisfying  the ordinary differential equation
\[
\sigma(D)\biggl(x',0,\xi',-i\frac{d}{dt}\biggr)u(t)=0,\quad
(x',\xi')\in S^*\partial M,
\]
that remain bounded as $t\to +\infty$ (see \cite{Hor3}). By
$\widehat{L}_+(D)$ let us denote some subspace in $C^\infty
(\partial M, E{\big|^d_{\partial M}})$ that is defined by a
pseudodifferential projection $Q$ on $\partial M$ with symbol
projecting on the subbundle $L_+(D)\subset E\big|^d_{\partial M}$
along the complementary subbundle $L_-(D)$ corresponding to
solutions decreasing at $-\infty$. Projections $Q$ are called
\textit{Calderon projections} for $D$.

The following Fredholm criterion is valid for boundary value
problems~\eqref{R}.

\begin{Theo}\label{teo5}
Boundary value problem~\eqref{R} defines a Fredholm operator for
$s>d/2$
\[
\begin{pmatrix}
 D\\Bj
\end{pmatrix}
 \,:\, H^s (M,E)\longrightarrow
\begin{array}{c}
H^{s-d}(M,F) \\
\oplus \\
P\, H^{s-1/2} (\partial M)
\end{array}
\]
if and only if two conditions are satisfied:

1) the principal symbol of $B$ is invertible on $S^* \partial
M\setminus \pi^* S^* X$
\[
\sigma_{\partial M} (B)\,:\,\op{Im}\sigma_{\partial M} (Q)
\longrightarrow \op{Im} \sigma (P);
\]

2) the operator symbol of $B$ is invertible on $S^* X$:
\[
\sigma_X(B)\,:\,\op{Im}\sigma_X(Q)\longrightarrow \op{Im} P.
\]
Here $Q$ is a Calderon projection for $D$ and for $d>1$ we denote by
$\sigma_{\partial M}(B)$ and $\sigma_X(B)$ the symbols of
the tuple $(B_0,\ldots, B_{d-1})$.
\end{Theo}

\begin{Rem}
It should be noted that the conditions  of the theorem use only
the principal symbol $\sigma (D)$, symbol of $B$ and the
projection $P$.
\end{Rem}

\begin{proof}
By the results of \cite{ScSS18}, the boundary value problem has
the Fredholm property if and only if an operator on the boundary
\[
(B_0,B_1, \ldots,B_{d-1})|_{\op{Im}Q}\,:\,\op{Im} Q
\rightarrow \op{Im}P
\]
has the Fredholm property. We apply the Fredholm criterion stated
in Theorem \ref{t5} to this operator in subspaces.

This readily shows that conditions 1) and 2) of the present
theorem are necessary and sufficient for the Fredholm property to
be valid.

\end{proof}

\noindent {\bf Example 3.} (Elliptic boundary value problem for
the Hirzebruch operator). Let $M^{4k}$ be an oriented
$4k$-dimensional manifold with boundary. Suppose that the boundary
is a product
\[
\partial M =X^{ev}\times Y^{odd}.
\]
In a neighborhood of the boundary choose a metric corresponding to
the Cartesian product $[0,1)\times X\times Y$. Consider the
Hirzebruch operator
\[
D_M=d+\delta\,:\,\Lambda^+(M)\longrightarrow \Lambda^-(M)
\]
on $M$, where $d$ and $\delta$ are the exterior derivative and its
adjoint and $\Lambda^{\pm}(M)$ are the spaces of dual and
antiselfdual forms on $M$ (e.\,g., see \cite{Pal1}). Then near the
boundary $D_M$ can be represented in the form
\begin{equation}\label{AAA}
 \frac{\partial}{\partial t}+
    \begin{pmatrix}
D_Y & D^*_X \\
D_X & -D_Y
  \end{pmatrix}
  ,
\end{equation}
where $D_X$ denotes the Hirzebruch operator on $X$
\[
D_X\,:\,\Lambda^+(X)\longrightarrow \Lambda^-(X)
\]
and $D_Y$ denotes the odd analog of the Hirzebruch  operator on Y
(see \cite{APS1}):
\[
D_Y\,:\,\Lambda^*(Y)\longrightarrow \Lambda^*(Y),
\]
\[
D_Y =\tau(d_Y+\delta_Y), \quad \tau\bigl|_{\Lambda^p(Y)}=
i^{(dimY+1)/2+p(p-1)}*, \quad \tau^2=id.
\]
It is elliptic and self-adjoint. Denote by $\Pi_+$
\[
\Pi_+\,:\,\Lambda^*(Y)\longrightarrow \Lambda^*(Y)
\]
the nonnegative spectral projection of $D_Y$ and $\Pi_-$ denote
the complementary projection.

Consider the boundary value problem
\begin{equation}\label{**}
\left\{
\begin{array} {l}
 \bigg[\dfrac{\partial}{\partial t}+
\begin{pmatrix}
D_Y & D^*_X\\
D_X & -D_Y
\end{pmatrix}
\bigg]  \begin{pmatrix} u\\v
\end{pmatrix}
=
\begin{pmatrix} f_1\\f_2
\end{pmatrix},
\vspace{0.2cm}\\
\Pi_+u\big|_{\partial M}=g_1\in\op{Im}\Pi_+,
\vspace{2mm}\\
\Pi_-v\big|_{\partial M}=g_2\in\op{Im}\Pi_-.
\end{array}\right.
\end{equation}

\begin{Prop}
Boundary value problem~\eqref{**} has the Fredholm property.
\end{Prop}

\begin{proof}
 1) According to Theorem \ref{teo5}, it suffices to check the
  invertibility of the corresponding principal and operator
  symbols.

 2) An elementary computation shows that the principal symbol of
 the Calderon projection $Q$ on $S^*\partial M$ is the  matrix
\[
\sigma_{\partial M}(Q)=\frac 12
\begin{pmatrix} 1+\sigma(D_Y) & \sigma^*(D_X)\\
\sigma (D_X) & 1-\sigma (D_Y)
\end{pmatrix}.
  \]
Hence, the operator symbol of the Calderon projection is
\[
\sigma_X(Q)=\frac 12
\begin{pmatrix} 1 & \sigma^*(D_X)\\
\sigma (D_X) & 1
\end{pmatrix}.
\]
We need to prove that the two maps
\[
\begin{array}{l}
\op{Im}\sigma_{\partial M}(Q)\overset{\pi_+\oplus
\pi_-}{\longrightarrow}\op{Im}(\pi_+\oplus \pi_- ),
\vspace{0.15cm}\\
\op{Im}\sigma_X(Q)\overset{\Pi_+\oplus
\Pi_-}{\longrightarrow}\op{Im} \Pi_+\oplus \op{Im}\Pi_-
\end{array}
\]
are isomorphisms (here $\pi_\pm$ denote the principal symbols of
$\Pi_\pm$). To prove that the two maps are isomorphisms one can
compute their compositions with the maps in the opposite direction
\[
\begin{array}{l}
\op{Im}(\pi_+\oplus \pi_- )\overset{\sigma_{\partial
M}(Q)}{\longrightarrow}\op{Im}\sigma_{\partial M}(Q),
\vspace{0.15cm}\\
\op{Im}\Pi_+\oplus
\op{Im}\Pi_-\overset{\sigma_X(Q)}{\longrightarrow}\op{Im}\sigma_X(Q).
\end{array}
\]
An explicit computation shows that the compositions
\[
(\pi_+\oplus \pi_-)\sigma_{\partial M}(Q), \quad \sigma_{\partial
M}(Q)(\pi_+\oplus \pi_-), \quad (\Pi_+\oplus \Pi_-)\sigma_X(Q),
\quad \sigma_X(Q)(\pi_+\oplus \pi_-)
\]
have trivial kernels. Thus, the boundary value problem~\eqref{**}
satisfies the assumptions of Theorem \ref{teo5} and, consequently,
has the Fredholm property.

\end{proof}

\section{Topological obstruction}\label{s6}

There is an obstruction to define a Fredholm boundary value
problem for a given elliptic operator $D$ on $M$.

\begin{Theo}
Suppose that an elliptic differential operator $D$ on $M$ has a
Fredholm boundary value problem of the form~\eqref{R}. Then the
principal symbol $\sigma(D)$ at the boundary has the following
property
\begin{equation}\label{8}
  \pi_![\sigma(D)|_{\partial M}]=0,
\end{equation}
where
\[
[\sigma(D)|_{\partial M}]\in K(T^*\partial
M\times\mathbb{R})\simeq K^1(T^*\partial M)
\]
is the difference construction and
\[
\pi_!\,:\, K^1(T^* \partial M)\rightarrow K^1(T^*X)
\]
is the direct image mapping in $K$-theory under the projection
$\pi$.
\end{Theo}

\begin{proof}

1) Let $(D,B,P)$ be an elliptic boundary value problem
  for $D$. This means that there is a Fredholm operator
\begin{equation}\label{la2}
  (B_0,B_1,\ldots,B_{d-1}):\,\op{Im}Q\rightarrow \op{Im P,}
\end{equation}
where $Q$ denotes a Calderon projection as previously.

The existence of a Fredholm isomorphism~\eqref{la2} implies that
the symbols of the two subspaces $\op{Im}Q$ and $\op{Im}P$ are
homotopic.

More precisely, consider both $Q$ and $P$ as projections in the
direct sum
\[
C^\infty(\partial M, E\big|^d_{\partial M})\oplus C^\infty
(\partial M, G).
\]
Then the homotopy of the principal symbols is defined as
\begin{equation}\label{alpha}
  q_{\partial M,\varphi}=\sigma(Q)\cos^2\varphi+\sigma(P)
  \sin^2\varphi+2\sigma(P)
  \sigma_{\partial M}(B)\sigma(Q)\sin\varphi\cos\varphi
\end{equation}
and a similar formula is valid for the operator symbol
\begin{equation}\label{beta}
 q_{X,\varphi}=\sigma_X(Q)\cos^2\varphi+P\sin^2\varphi+2P
 \sigma_X(B)\sigma_X(Q)\sin\varphi\cos\varphi.
\end{equation}
This is a homotopy of compatible symbols and we obviously have
at $\varphi=0$
$$
q_{\partial M,0}=\sigma(Q),\quad q_{X,0}=\sigma_X(Q)
$$
and for $\varphi=\pi/2$
$$
q_{\partial M,\pi/2}=\sigma(P),\quad q_{X,\pi/2}=P.
$$

2) On the other hand, let us represent the topological
  invariant in~\eqref{8} in analytic terms. It is well known that
  the element $[\sigma(D)|_{\partial M}]\in K^1(T^*\partial M)$
  can be expressed in terms of the symbol of the
  Calderon projection for $D$:
\begin{equation}\label{o}
  [\sigma(D)|_{\partial
  M}]=\big[(2\sigma(Q)(x,y,\xi,\eta)-1)\sqrt{\xi^2+\eta^2}+i\tau\big]
  \in K(T^*\partial  M\times\mathbb{R}),
\end{equation}
where the coordinates correspond to the fibration
$\pi\,:\,\partial M\rightarrow X$ and $\tau\in\mathbb{R}$ denotes
an additional variable. The element is understood in the sense of
the difference construction, since it is an isomorphism outside a
compact set.

Let us represent the element in Eq.~\eqref{o} as the difference
construction in the sense of \cite{AtSi4} for a family of elliptic
operators in the fibers, where the parameter space is the product
$T^*X\times\mathbb{R}$. More precisely, we will define a family of
elliptic operators parametrized by $B(T^*X\times\mathbb{R})$,
where $B$ denotes the unit ball bundle. The family will turn out
to be invertible over the spherical bundle
$S(T^*X\times\mathbb{R})$. Therefore, its analytic index is in the
following relative group
\[
K(B(T^*X\times\mathbb{R}),\,\,S(T^*X\times\mathbb{R}))\simeq
K(T^*X\times\mathbb{R}).
\]

The desired family of elliptic operators denoted by
$D(x,\xi,\tau)$ for the parameters $(x,\xi,\tau)\in
B(T^*X\times\mathbb{R})$ corresponds to a family of symbols
\begin{equation}\label{T}
  \bigg[2\sigma(Q)\left(x,y,\xi,\frac{\eta}{|\eta|}\bigl(1-\xi^2-\tau^2\bigr)
  \right)-1\bigg]\biggl(|\xi|+1-(\xi^2+\tau^2)\biggr)+i\tau.
\end{equation}
The invertibility of this symbol for $\xi^2+\tau^2\leq 1$ can be
verified by an explicit computation. It is clear that on the unit
spheres for $\xi^2+\tau^2=1$ this symbol is a vector bundle
isomorphism independent of $\eta$:
\[
\bigl(2\sigma(Q)(x,y,\xi,0)-1\bigr)|\xi|+i\tau.
\]

Hence, the difference construction for the family of elliptic
symbols~\eqref{T} is an element of the group $K(T^*\partial
M\times\mathbb{R})$.

One can show that this difference construction for the family
$D(x,\xi,\tau)$ coincides  with element~\eqref{o}. Hence, by the
Atiyah--Singer formula for families we obtain the desired
expression in analytic terms:
\begin{equation}\label{alo}
  \pi_!\big[\sigma(D)\big|_{\partial M}\big]=\op{ind} \big[D(\cdot,\cdot,\cdot)
  \big]\in
  K^1(T^*X),
\end{equation}
where
\begin{equation}\label{ali}
D(x,\xi,\tau)=
\bigg[2\sigma(Q)\left(x,y,\xi,\widehat{\frac{\eta}{|\eta|}}\bigl(1-
\xi^2-\tau^2\bigr)\right)-1\bigg]\bigl(|\xi|+1-({\xi^2+\tau^2})\bigr)+i\tau
\end{equation}
and the hat means that we have a family of operators in the
fibers.

We will show that this family of elliptic operators is homotopic
to an invertible family. Therefore, the index is zero in this
case.

Denote by $Q(\xi,\tau)$ the following family of pseudodifferential
operators in the fibers
\[
Q(\xi,\tau)=\sigma(Q)\left(x,y,\xi,\widehat{\frac{\eta}{|\eta|}}
\bigl(1-\xi^2-\tau^2\bigr)\right).
\]
On the sphere $\xi^2+\tau^2=1$ for $\xi\neq 0$ this is a family of
projections, while for other values of the parameters
$Q(\xi,\tau)$ is only an almost-projection. One can verify that
this property implies that~\eqref{ali} is elliptic on
$B(T^*X\times\mathbb{R})$ and invertible on
$S(T^*X\times\mathbb{R})$. Let us define the homotopy
$D_{\varphi}(x,\xi,\tau)$ by changing this family $Q(\xi,\tau)$
only:
\begin{equation}\label{S}
  D_{\varphi}(x,\xi,\tau)=\bigl(2Q_\varphi
  (\xi,\tau)-1\bigr)\bigl(|\xi|+1-({\xi^2+\tau^2})\bigr)+i\tau,
\end{equation}
such that $Q_\varphi (\xi,\tau)$ satisfies the above mentioned
property. The homotopy
$Q_\varphi(\xi,\tau),\,\,\varphi\in[0,\pi/2]$ is defined on the
spheres as
\[
Q_\varphi(\xi,\tau)=q_{X,\varphi}\biggl(\frac{\xi}{|\xi|}\biggr),\qquad
\text{for}\quad \xi^2+\tau^2=1,
\]
(the family $q_{X,\varphi}$ was defined in~\eqref{beta}), while
inside the balls for $\xi^2+\tau^2<1$ it suffices to define a
homotopy of the corresponding symbols:
\[
\sigma(Q_\varphi)(\xi,\tau)=q_{\partial
M,\varphi}\bigl(\xi,\eta(1-\xi^2-\tau^2)\bigr),
\]
(the principal symbols $q_{\partial M,\varphi}$ were defined
in~\eqref{alpha}). For brevity we omit the geometric variables $x,y$ in
the formulas.

One can verify that for this choice of $Q_\varphi$ the
operator  in~\eqref{S} is elliptic for $\xi^2+\tau^2\leq 1$ and is
an isomorphism for $\xi^2+\tau^2 =1$.

At the end of the homotopy for $\varphi=\pi/2$ the operator
$Q_\varphi(\xi,\tau)$ is a pseudodifferential projection that does
not depend on $\xi$ and $\tau$. Hence, for $\varphi=\pi/2$ the
family $D_\varphi (x,\xi,\tau)$ is invertible on the entire space
$B(T^*X\times\Bbb{R})$. Together
with~\eqref{alo} this yields the desired:
\[
\pi_!\big[\sigma(D)\big|_{\partial M}\big]=\big[\op{ind} D_{\frac
\pi 2}(\cdot,\cdot,\cdot)\big]=0.
\]

\end{proof}

\noindent {\bf Example 4.}
 Similar to Example 3, suppose that $\partial M^{4k}=X^{odd}
 \times Y^{ev}$ and consider the projection
 \[
X^{odd}\times Y^{ev}\overset{\pi}{\longrightarrow} X^{odd}
 \]
with an even-dimensional fiber $Y^{ev}$.

Now the Hirzebruch operator $D_M$ acquires the form
\[
\frac{\partial}{\partial t}+
\begin{pmatrix}
D_X & D^*_Y \\
D_Y & -D_X
\end{pmatrix}.
\]
 In contrast with the previous Example 3, the obstruction does not
 vanish in this case.

\begin{Lem}

\[
\pi_![\sigma(D)\big|_{\partial M}]=\op{ind}D_Y[\sigma(D_X)]\in
K^1(T^*X).
\]
\end{Lem}

\begin{proof}
This follows from the index formula for families.

\end{proof}

According to this Lemma, the Hirzebruch operator has no Fredholm
boundary value problems if $\op{ind} D_Y\neq 0$. When
$\op{ind}D_Y=0$, consider a boundary value problem
\begin{equation}\label{omega}
\left\{
\begin{array} {l}
 \bigg[\dfrac{\partial}{\partial t}+
\begin{pmatrix}
D_X & D_Y^*\\
D_Y & -D_X
\end{pmatrix}
\bigg]  \begin{pmatrix} u\\v
\end{pmatrix}
=
\begin{pmatrix} f_1\\f_2
\end{pmatrix},
\vspace{0.2cm}\\
u\big|_{\partial M}+D^*_Y(\triangle_Y+1)^{-1}v\big|_{\partial
M}=g.
\end{array}\right.
\end{equation}

\begin{Prop}
Boundary value problem~\eqref{omega} has the Fredholm property.
\end{Prop}

\begin{proof}
The check of the ellipticity of the boundary condition is similar
to the one in Example 3 and is left to the reader.

\end{proof}

\section{Appendix}

In this appendix we prove two important technical results.
First, we establish the boundedness in the Sobolev spaces
for a class of operators with operator-valued symbols.
Second, we prove the composition formula from Section \ref{s3}.

\begin{Theo}
\label{bounded}
Suppose that an operator-valued symbol $\widehat{p}\left( x,\xi \right) \in
\Psi ^d\left( Y\right) $ of order $d$ is defined for $x\in \Bbb{R}^n$,
vanishes outside a compact set and satisfies the estimates
\begin{equation}
\label{beto}
\left\| \left( 1+\Delta _X\right) ^N\left( \Delta _Y+\xi ^2\right) ^{\left(
s-d\right) /2}\widehat{p}\left( x,\xi \right) \left( \Delta _Y+\xi ^2\right)
^{-s/2}\right\| _{L^2\left( Y\right) \rightarrow L^2\left( Y\right) }\leq C_p
\end{equation}
uniformly in $x$ and $\xi $ for some $N>\left( n+\left| s-d\right| \right)
/2.$ Then the operator
$$
P:C^\infty \left(\Bbb{R}^n \times Y\right) \rightarrow
C^\infty \left(\Bbb{R}^n\times Y\right)
$$
with operator-valued symbol $\widehat{p}$
extends to an operator of order $d$ in
the Sobolev spaces and an estimate of its norm is valid
\[
\left\| P\right\| _{H^s\left(\Bbb{R}^n \times Y\right)
\rightarrow H^{s-d}\left(
\Bbb{R}^n\times Y\right) }\leq C_p\cdot C\left( s,d\right) ,
\]
where the constant $C\left( s,d\right) $ does not depend on the operator.
\end{Theo}

\begin{proof}
Let us estimate the norm $\left\| Pu\right\| _{H^{s-d}\left(
\Bbb{R}^n\times Y\right) }.$ In terms of the Fourier transform in $x$ this norm can
be represented as
\[
\left\| Pu\right\| _{H^{s-d}\left(\Bbb{R}^n \times Y\right) }^2=\left\| \left(
\Delta _Y+\zeta ^2\right) ^{\left( s-d\right) /2}\widehat{Pu}\left( \zeta
\right) \right\| _{L^2\left( \Bbb{R}^n,L^2(Y) \right) }^2.
\]
Since
\[
Pu=\int e^{ix\xi }\widehat{p}\left( x,\xi \right) \widehat{u}\left( \xi
\right) d\xi ,
\]
we obtain that
\[
\widehat{Pu}\left( \zeta \right) =\int \widetilde{\widehat{p}}\left( \zeta
-\xi ,\xi \right) \widehat{u}\left( \xi \right) d\xi ,
\]
where $\widetilde{\widehat{p}}$ denotes the Fourier transform in $x$ of the
operator symbol. Using this expression, we obtain
\[
\left\| Pu\right\| _{H^{s-d}\left(\Bbb{R}^n \times Y\right) }^2
=\left\| \int \left(
\Delta _Y+\zeta ^2\right) ^{\left( s-d\right) /2}\widetilde{\widehat{p}}%
\left( \zeta -\xi ,\xi \right) \widehat{u}\left( \xi \right) d\xi \right\|
_{L^2\left( \Bbb{R}^n,L^2\left( Y\right) \right) }^2.
\]
This shows that
$\left\| Pu\right\| _{H^{s-d}\left(\Bbb{R}^n \times Y\right) }^2$ is
equal to
\begin{equation}
\left\| \int \left(
\Delta _Y+\zeta ^2\right) ^{\left( s-d\right) /2}\widetilde{\widehat{p}}%
\left( \zeta -\xi ,\xi \right) \left( \Delta _Y+\xi ^2\right) ^{-s/2}\left(
\Delta _Y+\xi ^2\right) ^{s/2}\widehat{u}\left( \xi \right) d\xi \right\|
_{L^2\left( \Bbb{R}^n,L^2\left( Y\right) \right) }^2.  \label{expr}
\end{equation}
Let us now estimate the norm of the product in this formula. Consider
first the term
\[
A(\xi,\zeta)
=\left\| \left( \frac{\Delta _Y+\zeta ^2}{\Delta _Y+\xi ^2}\right) ^{\left(
s-d\right) /2}\left( \Delta _Y+\xi ^2\right) ^{\left( s-d\right) /2}%
\widetilde{\widehat{p}}\left( \zeta -\xi ,\xi \right) \left( \Delta _Y+\xi
^2\right) ^{-s/2}\right\| _{L^2\left( Y\right)\to L^2(Y) }.
\]
The estimate (\ref{beto}) of the theorem implies that
\[
A(\xi,\zeta)
\leq Const\left\| \left( \frac{\Delta _Y+\zeta ^2}{\Delta _Y+\xi ^2}\right)
^{\left( s-d\right) /2}\right\| _{L^2\left( Y\right)
\to L^2(Y)}\cdot C_p\frac 1{\left(
1+\left| \zeta -\xi \right| ^2\right) ^N}.
\]
On the other hand, the first term can be estimated by a Peetre type inequality
\[
\left\| \left( \frac{\Delta _Y+\zeta ^2}{\Delta _Y+\xi ^2}\right) ^{\left(
s-d\right) /2}\right\| _{L^2\left( Y\right)\to L^2(Y) }\leq
Const\left( 1+\left| \zeta
-\xi \right| ^2\right) ^{\left| s-d\right| /2}.
\]
Thus, the term $A(\xi,\zeta)$ is estimated as
\[
A(\xi,\zeta)\leq Const\cdot C_p
\frac 1{\left( 1+\left| \zeta -\xi \right| ^2\right)
^{N-\left| s-d\right| /2}}.
\]
We are now in a position to estimate (\ref{expr})
\[
\left\| Pu\right\| _{H^{s-d}
\left(\Bbb{R}^n \times Y\right)}^2\leq \left\| \int
A\left( \xi ,\zeta \right) \left\| \left( \Delta _Y+\xi ^2\right) ^{s/2}%
\widehat{u}\left( \xi \right) \right\| _{L^2\left( Y\right) }d\xi \right\|
_{L^2\left( \Bbb{R}^n\right) }^2.
\]
This gives
\[
\left\| Pu\right\| _{H^{s-d}\left( \Bbb{R}^n\times Y\right) }^2\leq Const\cdot
C_p\left\| u\right\| _{H^s\left( \Bbb{R}^n\times Y\right) }^2\cdot \int \frac{%
d\zeta }{\left( 1+\left| \zeta \right| ^2\right) ^{N-\left| s-d\right| /2}},
\]
where, for $N>\left( n+\left| s-d\right| \right) /2$ the last integral
converges. This proves the
boundedness of $P$ and the corresponding norm estimate.

\end{proof}

{\bf Proof of the composition formula.} The rest of the appendix contains the
proof of the composition formula
\[
\widehat{\sigma}_1\widehat{\sigma}_2\equiv
\widehat{\sigma_1\sigma_2}
\]
for two symbols $\sigma_1,\sigma_2\in\Sigma(M,\pi)$. The
comparison is valid modulo a compact operator in the scale of
Sobolev spaces.

\begin{Lem}\label{l1}
If the composition formula is valid for all
$\sigma_1\in\Sigma_0,\,\,\sigma_2\in\Sigma (M,\pi)$. Then it is
true in general.
\end{Lem}

\begin{proof}
Consider a pair of symbols $\sigma_1,\sigma_2\in\Sigma (M,\pi)$.
Choose an approximation of $\sigma_1$ by an element
$\sigma_\varepsilon$ of the dense subalgebra $\Sigma_0\subset
\Sigma (M,\pi)$ (see Proposition 1). Then the difference of the
corresponding operators is denoted by
\[
A_\varepsilon =\widehat{\sigma}_1-\widehat{\sigma}_\varepsilon
\]
and has norm less than $\varepsilon$ in any given Sobolev space,
provided the approximation is chosen appropriately. Thus, we
obtain
\[
\widehat{\sigma}_1\widehat{\sigma}_2=(\widehat{\sigma}_\varepsilon
+A_\varepsilon)\widehat{\sigma}_2.
\]
By the assumption of the Lemma we have
\[
(\widehat{\sigma}_\varepsilon
+A_\varepsilon)\widehat{\sigma}_2=\widehat{\sigma_\varepsilon\sigma_2}+A_
\varepsilon\widehat{\sigma}_2+K_\varepsilon
\]
(here $K_\varepsilon$ is a compact operator). Hence,
\[
\widehat{\sigma}_1\widehat{\sigma}_2=\widehat{\sigma_1\sigma_2}+A_
\varepsilon\widehat{\sigma}_2 +\widehat{\sigma}'+K_\varepsilon,
\]
where $\sigma'$ denotes the symbol
$(\sigma_\varepsilon-\sigma_1)\sigma_2$ with a small norm. As we
let $\varepsilon\to 0$, this equality shows that the difference
\[
\widehat{\sigma_1\sigma_2}-\widehat{\sigma}_1\widehat{\sigma}_2
\]
is a norm limit of a family of compact operators $K_\varepsilon$.
Therefore, the difference is compact as well.

\end{proof}

\begin{Lem}\label{l2}
The composition formula is valid for all symbols
$\sigma_1\in\Sigma_0,\,\,\sigma_2\in\Sigma(M,\pi)$.
Moreover, the error term has order $-1$ in the Sobolev spaces.
\end{Lem}
\begin{proof}
Obviously, it is sufficient to prove the formula when $\sigma_1$
either corresponds to an operator with a smooth symbol or a family
of pseudodifferential operators in the fibers. We consider the two
possibilities separately.

1) Let
\[
\sigma_1=(a(x,y,\xi,\eta),a(x,y,\xi,0))
\]
correspond to a smooth symbol $a(x,y,\xi,\eta)$ and
$\sigma_2=(p_M,p_X)$ be a general symbol. The composition of the
corresponding operators has the form
\[
\widehat{a}\bigl(\widehat{p}+\widehat{p_X-\widehat{\widetilde{p}}}_M)\bigr).
\]
The desired composition formula will be proved if we show that
\begin{equation}\label{i}
  \widehat{a}\widehat{p}_M\equiv\widehat{ap}_M
  \end{equation}
and for a compact operator-valued symbol
$q_X=p_X-\widehat{\widetilde{p}}_M$ a similar comparison is valid
\begin{equation}\label{j}
\widehat{a}\widehat{q}_X\equiv\widehat{\widetilde{a}q_X}.
\end{equation}

a) To prove~\eqref{i} we have to show that the usual composition
formula for $\psi DO$'s on $M$ remains valid in this case (since
the expressions for the operators in~\eqref{i} do not contain
additional  operator-valued components).

The usual proof of the composition formula (e.g., see
\cite{EgSc1}, p. 32) can be repeated verbatim to establish~\eqref{i}.
Moreover, the composition formula is valid modulo an operator of
order $-1$ in the scale of Sobolev spaces.

b) To prove~\eqref{j} let us represent operator $\widehat{a}$ as
an operator on $X$ with an operator-valued symbol (see
\cite{Luk1}):
\[
a'(x,\xi)=a\left(x,y,\xi,-i\frac{\partial}{\partial
y}\right)\,:\,L^2(Y_x)\longrightarrow L^2(Y_x).
\]
Let us first prove that the composition  formula is satisfied
modulo a compact operator in $L^2(M)$.
From the composition formula for operators with operator-valued symbols (see
\cite{Luk1}) we obtain
\begin{equation}
\label{trio}
\widehat{a}\widehat{q}_X-\widehat{a'q_X}\equiv0.
\end{equation}
On the other hand, the right-hand side of the desired formula
(\ref{j}) in this case is equal to $\widehat{\widetilde{a}\, q_X}$.

To prove that the two expressions $\widehat{a'q_X}$ and
$\widehat{\widetilde{a}\,q_X}$ differ by a compact operator, it
suffices to show (see \cite{Luk1}) that the $L^2$-norm of the compact
operator-valued symbol
\begin{equation}\label{V}
\bigg[a\left(x,y,\xi,-i\frac{\partial}{\partial y
}\right)-a\bigl(x,y,\xi,0\bigr)\bigg]q_X(x,\xi)
\end{equation}
tends to zero as $|\xi|\to \infty$. Indeed, the symbol of the
operator in the square brackets is equal to
$a(x,y,\xi,\eta)-a(x,y,\xi,0)$ and can be estimated as:
\[
|a(x,y,\xi,\eta)-a(x,y,\xi,0)|\leq C\frac{|\eta|}{|\xi|+|\eta|+1}.
\]
A similar estimate is valid for the derivatives  in $y$ (with a
possibly different constant $C$).

From this estimate it follows that the norm of the corresponding
operator is estimated as follows
\[
\bigg\|a\left(x,y,\xi,-i\frac{\partial}{\partial
y}\right)-a\bigl(x,y,\xi,0\bigr)\bigg\|_{H^1(Y_x)\rightarrow
H^0(Y_x)}\leq \frac{C'}{1+|\xi|}.
\]

Hence, the norm of the symbol~\eqref{V} in $L^2$ is bounded by:
\begin{multline*}
\bigg\|\bigg[a\left(x,y,\xi,-i\frac{\partial}{\partial
y}\right)-a\bigl(x,y,\xi,0\bigr)\bigg]\cdot
q_X(\xi)\bigg\|_{L^2(Y_x)\rightarrow L^2(Y_x)}\leq \\
\|q_X\|_{L^2(Y_x)\rightarrow H^1(Y_x)}
\times\bigg\|a\left(x,y,\xi,-i\frac{\partial}{\partial
y}\right)-a\bigl(x,y,\xi,0\bigr)\bigg\|_{H^1(Y_x)\rightarrow
L^2(Y_x)}\leq \frac{C''}{1+|\xi|}.
\end{multline*}
Hence, it tends to zero as desired. This proves that the
composition formula is valid in this case modulo a compact operator
in $L^2(M)$. However, we need to prove a stronger statement
that the error term has order $-1$ in the {\em Sobolev scale}.
To prove this, it suffices to show that the difference in (\ref{trio})
and the operator with symbol (\ref{V}) have order $-1$.

First, similar to the proof of Theorem \ref{t1}, one can obtain
for the symbol in (\ref{V}) the estimate of the form
$$
\left\|(\sqrt{\Delta_Y}+|\xi|)^{s+1}
(a'-\widetilde{a})q_X (\sqrt{\Delta_Y}+|\xi|)^{-s}
\right\|_{L^2(Y_x)\to L^2(Y_x)}\le C
$$
and a similar estimate for its derivatives in $x$. Thus, similar to the
Theorem \ref{bounded}, this proves that the symbol
$(a'-\widetilde{a})q_x$ gives an operator of order $-1$ in the Sobolev scale.

Second, to estimate the difference (\ref{trio}) one should estimate the corresponding
error term in the composition formula (e.g., see \cite{EgSc1}, p. 32).
This can also be done. The details are left to the reader.

2) Let us now verify the composition formula for a symbol
\[
\sigma_1=(b(x,y,\eta),B(x))
\]
of a family of operators $B(X)$ in the fibers and
$\sigma_2=(p_M,p_X)$ as before. In this case the proof of the
composition formula amounts to verifying the two comparisons (here
$q_X=p_X-\widehat{\widetilde{p}}_M$ as before)
\begin{equation}\label{i'}
B\widehat{q_X}\equiv\widehat{Bq_X},
\end{equation}
\begin{equation}\label{j'}
B\widehat{p_M}\equiv\widehat{bp_M}+\widehat{B\widehat{\widetilde{p}}_M-
\widehat{b\widetilde{p}}_M}.
\end{equation}

Concerning the first composition, it is easy to see that it is
satisfied exactly.

Let us establish comparison~\eqref{j'}. To this end, we rewrite
the terms $B\widehat{p}_M$ and $\widehat{bp}_M$ as operators on
the base $X$ with operator-valued symbols. We introduce the
following notation. For a symbol $\sigma(x,y,\xi,\eta)$ denote the
operator-valued symbol $\sigma(x,y,\xi,-i{\partial}/{\partial y})$
by $\sigma'$. Then~\eqref{j'} can be rewritten as
\[
\widehat{Bp'_M}\equiv\widehat{(bp_M)'}+\widehat{q}_X,
\]
where the operator symbol $q_X$ is equal to
\[
q_X=B\widehat{\widetilde{p}}_M-\widehat{b\widetilde{p}}_M.
\]

Thus, to prove that the comparison~\eqref{j'}  is  valid
modulo a compact operator, it suffices to show that
the operator-valued symbol
\begin{equation}\label{J}
Bp'_M-(bp_M)'-B\widehat{\widetilde{p}}_M+\widehat{b\widetilde{p}}_M
\end{equation}
is a compact operator for all $(x,\xi)\in T^*X$ and its norm tends
to zero as $|\xi|\to \infty$.

First of all, the compactness of~\eqref{J} is easy to obtain
taking the symbols of the operators involved:
\begin{equation}
 \sigma(Bp'_M-(bp_M)'-B\widehat{\widetilde{p}}_M+\widehat{b\widetilde{p}_M})
= 0.
\end{equation}
Here we omitted for brevity the variables $x,y$ in the symbols.

Let us now show that the norm of the symbol~\eqref{J} tends to
zero as $|\xi|\to \infty$. To this end, we rewrite this symbol in
the form
\begin{equation}\label{J'}
  B(p_M-\widetilde{p}_M)'-(b(p_M-\widetilde{p}_M))'=(B-b')(p_M-
  \widetilde{p}_M)'+b'(p_M-\widetilde{p}_M)'-(b(p_M-\widetilde{p}_M))'.
\end{equation}
For the difference $p_M-\widetilde{p}_M$ the following estimate
can be obtained
\[
|p_M(\xi,\eta)-\widetilde{p}_M(\xi,\eta)|\leq
C\frac{|\eta|}{|\xi|+|\eta|+1}.
\]

Similar to the previous part of the proof this shows that the
operator-valued symbol $(B-b')(p_M-\widetilde{p}_M)'$
in~\eqref{J'} tends to zero for $\xi\to \infty$. Let us estimate
the remaining term
\[
b'(p_M-\widetilde{p}_M)'-(b(p_M-\widetilde{p}_M))'
\]
by means of the usual composition formula on the fiber. It is a
pseudodifferential operator of order $-1$ with a symbol
$C(\xi,\eta)$ estimated as
\[
|C(\xi,\eta)|\leq\frac{C_1}{|\xi|+|\eta|+1}.
\]
Thus, the norm of the corresponding operator-valued symbol tends
to zero as $|\xi|\to \infty$.  The proof that the composition formula
is valid modulo an operator of order $-1$ can be obtained similar to the
previous part of the proof.

This completes the proof of the composition formula.

\end{proof}


\end{document}